\documentstyle{amsppt}

\input psfig

\topmatter
\title
On the frequency of vanishing of quadratic twists of modular L-functions
\endtitle
\author
J.B.~Conrey \\
J.P.~Keating\\
M.O.~Rubinstein\\
N.C.~Snaith
\endauthor
 \address                                                  
 American Institute of Mathematics,
Palo Alto, CA 94306
\endaddress
\address
Basic Research Institute in the Mathematical Sciences, Hewlett-Packard
Laboratories, Bristol BS34 8QZ U.K.\endaddress
\address University of Texas, Austin, TX\endaddress
\address University of Bristol, Bristol U.K. \endaddress

\abstract
We present theoretical and numerical evidence for a random matrix
theoretical approach to a conjecture about vanishings of quadratic twists
of certain L-functions.  \endabstract
\endtopmatter
\document

\NoBlackBoxes           

\documentstyle{amsppt}
\document

In this paper we 
\footnote{Supported in part by a Focused Research Group grant from the NSF.}
present some evidence that methods from random matrix
theory can give insight into the frequency  of vanishing for quadratic
twists of modular $L$-functions. The central question   is the
following: given a holomorphic newform $f$ with integral coefficients and
associated $L$-function $L_f(s)$,  for how many fundamental
discriminants $d$ with
$|d|\le x$, does $L_f(s,\chi_d)$, the $L$-function twisted by the real,
primitive, Dirichlet character associated with the discriminant $d$,
vanish at the center of the critical strip to order at least 2?

This question is of particular interest in the case that the $L$-function
is associated with an elliptic curve, in light of the conjecture of Birch
and Swinnerton-Dyer. This case corresponds to weight $k=2$.
We will focus on this case for most of the paper, though 
we do make some remarks about higher weights (see (26) and below).

Suppose that $E/Q$ is an elliptic curve with associated $L$-function 
$$L_E(s)=\sum_{n=1}^\infty \frac{a^*_n}{n^s} \tag1$$
for $\Re s>1$.  
Then, as a consequence of the Taniyama-Shimura conjecture, recently solved
by Wiles, Taylor, ([W], [TW]), and  Breuil, Conrad, and Diamond, $L_E$ is entire
and satisfies  a functional equation
$$\bigg(\frac{2\pi}{\sqrt{N}}\bigg)^{-s}\Gamma(s)L_E(s)=\Phi_E(s)=w_E\Phi(1-s) \tag2$$
where $N$ is the conductor of $E$ and $w_E=\pm1$ is called the sign of the
functional equation.  Note that we have normalized the coefficients
$a_n^*$ so that the functional equation of $L_E$ relates the values at
$s$ and $1-s$. The numbers $a_n^*$ satisfy $|a_n^*|\le d(n)$, where
$d(n)$ is the number of positive divisors of $n$, and 
$a_n=\sqrt{n}a_n^*$ is an integer (related to the numbers of points on
$E$ mod $p$
for primes $p$ which divide $n$).  
Let $d$ represent a fundamental discriminant 
and let $\chi_d$ be the associated
quadratic character (i.e. $\chi_d(n)=\left(\frac d n \right)$, the
Kronecker symbol). 
We assume, for simplicity, that $(d,N)=1$.
Then the twisted $L$-function is
$$L_E(s,\chi_d)=\sum_{n=1}^\infty \frac {a_n^*\chi_d(n)}{n^s}.\tag3 $$
This $L$-function has functional equation 
$$\bigg(\frac{2\pi}{|d|\sqrt{N}}\bigg)^{-s}\Gamma(s)L_E(s,\chi_d)=\Phi_E(s,\chi_d)=w_E\chi_d(-N)\Phi(1-s,\chi_d);\tag4$$
it is actually the $L$-function of another elliptic curve, namely the
quadratic twist of $E$ by $d$.

The conjecture of Birch and
Swinnerton-Dyer predicts that the order of vanishing of the $L$-function
of
an elliptic curve at the central critical point $s=1/2$ is the same as the
rank of the Mordell-Weil group of the elliptic curve.  Thus it
would be of interest to 
find an asymptotic formula for 
$$V_E(x):=\sum\Sb |d|\le x\\ w_E\chi_d(-N)=1  \\ L_E(1/2,\chi_d)=0\endSb
1\tag5$$
since this would potentially be counting how often the twists of a given
elliptic curve have rank at least 2.  
(Note that we have restricted the sum to twists for which the sign of the
functional equation is +1; these $L$-functions vanish to order at least 2
because of the symmetry implied by the functional equation.)
Goldfeld [G] has predicted that $V_E(x)=o(x)$.  More specifically, he predicts
that  asymptotically, 1/2 of all twists will have rank 0 and 1/2 of all
twists will have rank 1; consequently  ranks 2 and higher should be
infrequent.
  We will give a more precise conjecture about the frequency of twists
with ranks at least 2.

Sarnak has predicted that $V_E(x)$ should be about $x^{3/4}$.  His
reasoning has to do with the formulas of Waldspurger [Wa], Shimura [Sh], and Kohnen
and Zagier [KZ] which relate the value 
of $L_E(1/2,\chi_d)$ to the Fourier coefficient of a half-integral weight
modular form.  Roughly,
$$  L_E(1/2,\chi_d) = \kappa_E c_E(|d|)^2/\sqrt{d} \tag6$$
where $\kappa_E$ depends only on $E$ and where the integers $c_E(|d|)$ 
are the Fourier coefficients of a half-integral weight form.
The Ramanujan conjecture for these coefficients predicts that 
$c_E(|d|)\ll |d|^{1/4+\epsilon}$ for every $\epsilon>0$.  If $c_E(d)$
takes on each integer value up to $|d|^{1/4}$ about the same number of
times for $|d|\le x$, then it should take the value 0 about $x^{3/4}$
times.

Using random matrix theory, we would like to give a conjecture of the form  
$$V_E(x) \sim b_E x^{3/4} (\log x)^{e_E}\tag7$$
for certain constants $b_E$ and $e_E$.  
The basic idea is to regard the family 
$$ \Cal{F}_{E^+}=\{L_E(s,\chi_d):w_E\chi_d(-N)=+1\}\tag8$$
as an orthogonal family, in the sense of the families introduced by Katz
and Sarnak ([KS1], [KS2]).
More specifically, this family conjecturally has symmetry type
$O^+$.  Thus, for example, we believe that the statistics of the low lying
zeros of the $L$-functions 
in {\it this family} will match the statistics of eigenvalues near 1 of
the
matrices in
$SO(2N)$.  

The point of departure for our conjectures is the work of Keating and
Snaith [KeSn1] and [KeSn2] (see also [BH] and [CF]) which indicates that
the moments 
$$M_E(T,k)=\frac 1 {T^*}\sum\Sb |d|\le T\\L_E(s,\chi_d)\in \Cal{F}_{E^+}\endSb
L_E(1/2,\chi_d)^k$$
(with $T^* =\sum\Sb |d|\le T\\L_E(s,\chi_d)\in \Cal{F}_{E^+}\endSb 1$) 
apparently behave like the moments of the
characteristic polynomials of matrices 
in $SO(2N)$ where $N$ is of size $\log T$.  
Precisely, they   conjecture that
$$M_E(T,k) \sim g_k(O^+)a_k(E)(\log T)^{k(k-1)/2}\tag9$$
where 
$$
    g_k(O^+)=   2^{k(k+1)/2} \prod_{\ell=1}^{k-1}\frac{\ell!}{2\ell!}
$$
for integer $k$ and 
$$
\align
&a_k(E)= \tag10 \\
&\prod_p \left(1-\frac{1}{p}\right)^{k(k-1)/2}\left(\frac
{\left(1-\frac{a_p}{p}+\frac{1}{p}\right)^{-k}+\left(1+\frac{a_p}{p}+\frac{1}{p}\right)^{-k}}{2}
\frac{p}{p+1}+\frac{1}{p+1}\right).
\endalign
$$
This conjecture arises from arithmetical considerations together with the
fact from random matrix theory that the moments of the characteristic
polynomials of matrices in SO(2$N$), evaluated at the point 1, averaged
over the group can be explicitly evaluated. Thus,
$$
\align
M_O(N,s)&=\int_{SO(2N)}|\det(U-I)|^s~dU \tag11 \\
        &=2^{2Ns}\prod_{j=1}^N\frac{\Gamma(N+j-1)\Gamma(s+j-1/2)}{\Gamma(j-1/2)\Gamma(s+j+N-1)}
\endalign
$$
where $dU$ is the Haar measure for SO(2$N$).
The connection with $g_k$ is that
$$M_O(N,k)\sim g_k(O^+) N^{k(k-1)/2}\tag12$$
as $N\to \infty$.
Note that the formula for $g_k(O^+)$ can be extended to all real $k$ by
$$g_k(O^+)=2^{k^2/2}\frac{G(1+k) \sqrt{\Gamma(1+2k)}}{\sqrt{G(1+2k)\Gamma(1+k)}}\tag13$$
where $G$ is Barnes' Double Gamma function.

Continuing to follow Keating and Snaith ([KeSn2], equations (74) - (81)),
we 
observe that knowledge of all the complex moments of characteristic
polynomials,
evaluated at 1, in $SO(2N)$ gives complete information about the density
function for the distribution of values of
the characteristic polynomials at this point. Specifically, the latter
is the Mellin transform of the former:
$$P_O(N,x)=\frac{1}{2\pi i x}\int_{(c)}M_O(N,s)x^{-s} ds \tag14$$
where $(c)$ denotes the vertical line path from $c-i\infty$ to
$c+i\infty$; this formula is valid for all real $x$.  Note that
$P_O(N,x)~dx$ gives the probability that $\det(U-I)=x$ for an element $U$
of $SO(2N)$.

For small positive $x$, the pole of $M_O(N,s)$ at $s=-1/2$ determines the
dominating behavior of $P_O(N,x)$.  In fact, we see that
$$
\align
P_O(N,x) \sim &x^{-1/2} 2^{-N}\Gamma(N)^{-1}\prod_{j=1}^N
\frac{\Gamma(N+j-1)\Gamma(j)}{\Gamma(j-1/2)\Gamma(j+N-3/2)} \tag15 \\
&:=x^{-1/2} h(N)
\endalign
$$
as $x\to 0^+$. As $N\to \infty$, 
$$h(N)\sim  2^{-7/8}G(1/2)\pi^{-1/4} N^{3/8}.\tag16$$

An interpretation of the above is that the probability that an element of
SO(2$N$) has a characteristic polynomial whose value at 1 is $X$ or
smaller is
$$\sim \int_0^X x^{-1/2} h(N)~dx =2 X^{1/2}h(N).\tag17$$

We apply this reasoning to the values of $L_E(1/2,\chi_d)$. In particular,
by (6) the fact that the $c_E(|d|)$ are integers implies that these values
are discretized.  If, for example, it is known that
$$L_E(1/2,\chi_d)< \kappa_E/\sqrt{|d|}, \tag 18$$
 then it follows that $L_E(1/2,\chi_d)=0.$  Similarly, if
$$\kappa_E/\sqrt{|d|}\leq L_E(1/2,\chi_d)< 4 \kappa_E/\sqrt{|d|}, \tag
19$$
then it must be the case that 
$ L_E(1/2,\chi_d)= \kappa_E/\sqrt{|d|}.$

We assume now that the distribution of values of $L_E(1/2,\chi_d)$ will
behave like the values of 
the determinants of random orthogonal matrices with some suitable
restrictions and use this assumption to conjecture results about the
frequency of vanishing of $L_E(1/2,\chi_d)$.  The restrictions we have in
mind are of an arithmetical nature. First of all, we want to include the
arithmetical factor
$a_k(E)$.  We expect that
$$M_E(T,s)\sim a_s(E)M_O(N,s)\tag 20$$
with $N\sim \log T$.  Thus, using 
$$P_E(T,x)=\frac{1}{2\pi i x}\int_{(c)}M_E(T,s)x^{-s} ds \tag21$$
together with (20), we reiterate the conjecture of [KeSn2] (equation
(81)) that
$$
\align
P_E(T,x) &\sim a_{-1/2}(E) x^{-1/2} 2^{-N}\Gamma(N)^{-1}\prod_{j=1}^N
\frac{\Gamma(N+j-1)\Gamma(j)}{\Gamma(j-1/2)\Gamma(j+N-3/2)} \tag22 \\
&= a_{-1/2} x^{-1/2} h(N)
\endalign
$$
should approximate for small $x$ the probability density function for values of
$L_E(1/2,\chi_d).$
Of course, this formula cannot be too accurate as we have already remarked
that the values of $L_E(1/2,\chi_d)$ are discretized.
The precise nature of this discretization is somewhat involved; it
involves the constant $\kappa_E$ for which we have explicit formulas, but
it also involves the coefficients $c_E(|d|)$ whose arithmetic nature is
difficult to describe.  Simplistically, we would like to use 
the integral of (22), as in (17), and (18) to predict that
$$
\align
&\#\{|d|\le T: L_E(1/2,\chi_d)=0, L_E(s,\chi_d)\in \Cal{F}_{E^+} \} \tag23 \\
&\sim \frac{8}{3}\sqrt{\kappa_E}a_{-1/2} T^*/T^{1/4} h(N)
\endalign
$$  
with $N\sim\log T$. 
However, the $c(|d|)$ are divisible by some predetermined powers of 2
which change this discretization. For example, in the case of the
congruent number curve $E_{32}: y^2=x^3-x$, the number $c(|d|)$ is
divisible by $\tau(d)$ for squarefree $d$ where $\tau(d)$ is the number of
divisors of $d$. 
Thus, if  
$$L_E(1/2,\chi_d)< \kappa_{E_{32}}\tau(d)^2 /\sqrt{|d|}, \tag 24$$
then $L_{E_{32}}(1/2,\chi_d)=0$.
If we introduce the factor $\tau(d)$ into (23) it will raise the expected
frequency of vanishing by a factor of about $\log T$ giving a total order
of magnitude $T^{3/4}(\log T)^{11/8}$
for the frequency of vanishing. We expect this to be the correct order of
magnitude but are not able to say yet what constant we expect.

If we restricted to prime twists ($|d|=p$) then the extra powers of 2 are
not so significant. This leads to
\proclaim{Conjecture 1} Let $E$ be an elliptic curve defined over
$Q$.  Then
there is a constant $c_E>0$ such that  
$$\sum\Sb p\le T\\L_E(1/2,\chi_p)=0 \\ L_E(s,\chi_p)\in \Cal{F}_{E^+}\endSb 1 \sim c_E T^{3/4}(\log
T)^{-5/8}. \tag25$$
\endproclaim

We will return to a discussion of the value of $c_E$ in another paper.

We remark that our arguments apply equally well to newforms $f$ of weight
4 with integral Fourier coefficients and even functional equation.  
Here we have a different
discretization and expect that 
$$\sum\Sb p\le T\\L_f(1/2,\chi_p)=0 \endSb 1 
\sim c_f T^{1/4}(\log T)^{-5/8}. \tag26$$
Note in particular that the exponent on $T$ is now $1/4$. For a newform of
weight 6 or higher, we expect that there will be at most a finite number
of twists that vanish.
We have some numerical evidence to support these conjectures.

For the remainder of this  paper we would like to discuss a numerical
experiment which allows us to skirt the delicate issue of the arithmetic
nature of the $c(|d|)$. For a prime $p$ we consider the ratios
$$R_p(T)= \big(\sum \Sb |d|\le T \\
\chi_d(p)=1\\L_E(1/2,\chi_d)=0\\L_E(s,\chi_d)\in \Cal{F}_{E^+}\endSb1\big)
 /\big(\sum \Sb |d|\le T \\ \chi_d(p)=-1\\L_E(1/2,\chi_d)=0\\
L_E(s,\chi_d)\in \Cal{F}_{E^+}\endSb1
\big).\tag 27$$
By considering this ratio, the
powers of $T$, of $\log T$, and the constants intrinsic to the curve $E$
should all cancel out. 

More generally, let
$$ Q_p(k)=\lim_{T\to \infty}
\sum \Sb |d|\le T \\ \chi_d(p)=1
\\L_E(s,\chi_d)\in \Cal{F}_{E^+}\endSb L_E(1/2,\chi_d)^k/
\sum \Sb |d|\le T \\ \chi_d(p)=-1
\\L_E(s,\chi_d)\in \Cal{F}_{E^+}\endSb L_E(1/2,\chi_d)^k \tag
28$$
assuming that this limit exists. 
What is a reasonable conjecture for $Q_p(k)$?

Using standard techniques from analytic number theory (see [I]), we can
evaluate $Q_p(1)$.
Based on the analysis involved in such an evaluation, we expect that
$$Q_p(k)=
\frac{(p+1+a_p)^k}{(p+1-a_p)^k} \tag
29$$
where $a_p$ is the $p$-th Fourier coefficient of the modular form
associated with $E$. The heuristics are as follows. Consider
either sum that appears in (28)
$$
\align
\sum \Sb |d|\le T \\ \chi_d(p)=\pm 1
\\L_E(s,\chi_d)\in \Cal{F}_{E^+}\endSb L_E(1/2,\chi_d)^k
&=
\sum \Sb |d|\le T \\ \chi_d(p)=\pm 1
\\L_E(s,\chi_d)\in \Cal{F}_{E^+}\endSb
\bigg(\sum_{n=1}^\infty \frac {a_n\chi_d(n)}{n}\bigg)^k \tag30 \\
&=
\sum \Sb |d|\le T \\ \chi_d(p)=\pm 1
\\L_E(s,\chi_d)\in \Cal{F}_{E^+}\endSb
\sum_{n=1}^\infty \frac {b_n\chi_d(n)}{n}
\endalign
$$
where
$$
b_n = \sum_{n=n_1 n_2 \cdots n_k} a_{n_1} a_{n_2} \cdots a_{n_k},
$$
the sum being over all ways of writing $n$ as a product of $k$ factors.
Summing over $d$, the main contribution to (30) comes from those $n$'s that are
of the form $p^r m^2$, i.e. a power of $p$ times a perfect square
(since, unless $(d,n)>1$, these always have $\chi_d(n) = \chi_d(p)^r$ while, for 
other $n$'s, we get cancellation as we sum over $d$). So, the main 
contribution to (30) is roughly
$$
\sum \Sb |d|\le T \\ \chi_d(p)=\pm 1
\\L_E(s,\chi_d)\in \Cal{F}_{E^+}\endSb
\sum_{p^r m^2} \frac {b_{p^r m^2}\chi_d(p)^r}{p^r m^2}
$$
But, $b_{uv}=b_u b_v$ when $(u,v)=1$, so the inner sum above equals
$$
\sum_{(m,p)=1}
\frac {b_{m^2}}{m^2}
\sum_{r=0}^\infty \frac {b_{p^r}\chi_d^r(p)}{p^r}.
$$
Now,
$$
\sum_{r=0}^\infty \frac {b_{p^r}\chi_d^r(p)}{p^r}
=
\bigg(\sum_{r=0}^\infty \frac {a_{p^r}\chi_d^r(p)}{p^r}\bigg)^k.
\tag31
$$
Further, one can write an explicit formula for $a_{p^r}$ in
terms of $a_p$ and $p$. Assume that $E$ has good reduction mod $p$. Writing
the corresponding Euler factor
$(1-a_p p^{-s}+p^{1-2s}) = (1- \alpha p^{-s})(1-\beta p^{-s})$,
with $\alpha+\beta=a_p$ and $\alpha \beta = p$, we have, using
partial fractions, 
$a_{p^r}= (\alpha^{r+1}-\beta^{r+1})/(\alpha-\beta)$. Substituting this into
(31), and summing the geometric series we get
$$
\bigg(\sum_{r=0}^\infty \frac {a_{p^r}\chi_d^r(p)}{p^r}\bigg)^k
=
\frac {p^{2k}}{(p+1-\chi_d(p) a_p)^k}.
$$
Hence, the only difference in the numerator and denominator of
(28), asymptotically, is a factor of 
$$
\frac{(p+1+a_p)^k}{(p+1-a_p)^k}.
$$

Thus, by the random matrix theory considerations above
(as in (15)), 
taking $k=-1/2$ leads to
\proclaim{Conjecture 2}
With $R_p(T)$ defined as above, and $E$ having good reduction mod $p$,
$$R_p =\lim_{T\to \infty }R_p(T)= 
      \sqrt{\frac{p+1-a_p}{p+1+a_p}}.$$
\endproclaim

Note that the number $N_p$ of points on the elliptic curve $E$ over the
finite field $F_p$ of $p$ elements can be computed as 
$$ N_p=p+1-a_p$$
so that the above ratio is
the square-root of the ratio of the number of points on $E/F_p$ to the
number of points on $E^ \chi/F_p$ for any character $\chi$ with
$\chi(p)=-1$.
 
We conclude with some numerical evidence to support these conjectures. We
consider three elliptic curves which we call $E_{11}$, $E_{19}$, and
$E_{32}$. These are associated to the unique newforms of weight 2 and
levels 11, 19, and 32 respectively. The values of $L(s,\chi_d)$ were
evaluated by computing the $c(|d|)$'s of the corresponding weight $3/2$
forms. For the level 11 and 19 curves, we only computed these for $d<0$ and
$d$ odd. For the level 32 curve, we computed these for all odd $d$. The relevant form
for the level 32 case is described in [Ko]. The forms for
the level 11 and 19 cases can be computed according to [Gr] and were
given to us by Fernando Rodriguez-Villegas.

\topinsert
\hbox{\vbox{
\halign{\indent\hfil # 
&  \quad \hfil # & \hfil # & \quad \hfil #& \hfil #& \quad \hfil #& \hfil # \cr
$p$ & conjectured & data  & conjectured & data &conjectured & data \cr
    & $R_p$ for $E_{11}$& for $E_{11}$ & $R_p$ for $E_{19}$ & for $E_{19}$& $R_p$ for $E_{32}$& for $E_{32}$ \cr
 \cr
3 & 1.2909944 & 1.2774873 & 1.7320508 & 1.7018241 & 1 & 0.99925886 \cr
5 & 0.84515425 & 0.84938811 & 0.57735027 & 0.57825622 & 1.4142136 & 1.4113424 \cr
7 & 1.2909944 & 1.288618 & 1.1338934 & 1.134852 & 1 & 1.0003445 \cr
11 &  & 0 & 0.77459667 & 0.76491219 & 1 & 1.0001457 \cr
13 & 0.74535599 & 0.73266305 & 1.3416408 & 1.3632977 & 0.63245553 & 0.61626177 \cr
17 & 1.118034 & 1.1282072 & 1.183216 & 1.196637 & 0.89442719 & 0.88962298 \cr
19 & 1 & 1.000864 &  & 0 & 1 & 1.0006726 \cr
23 & 1.0425721 & 1.0470095 & 1 & 0.99857962 & 1 & 1.0000812 \cr
29 & 1 & 0.99769402 & 0.81649658 & 0.80174375 & 1.4142136 & 1.4615854 \cr
31 & 0.80064077 & 0.78332934 & 1.1338934 & 1.143379 & 1 & 1.0008405 \cr
37 & 0.92393644 & 0.91867671 & 0.9486833 & 0.94311279 & 1.0540926 & 1.0603105 \cr
41 & 1.2126781 & 1.2400086 & 1.1547005 & 1.1683113 & 0.78446454 & 0.76494748 \cr
43 & 1.1470787 & 1.1642671 & 1.0229915 & 1.0229106 & 1 & 1.0006774 \cr
47 & 0.84515425 & 0.82819492 & 1.0645813 & 1.0708874 & 1 & 0.99951502 \cr
53 & 1.118034 & 1.1332312 & 0.79772404 & 0.77715638 & 0.76696499 & 0.74137107 \cr
59 & 0.91986621 & 0.91329134 & 1.1055416 & 1.1196252 & 1 & 0.99969828 \cr
61 & 0.82199494 & 0.79865031 & 1.0162612 & 1.0199932 & 1.1766968 & 1.1996892 \cr
67 & 1.1088319 & 1.1216776 & 1.0606602 & 1.0705574 & 1 & 1.0002831 \cr
71 & 1.0425721 & 1.0497774 & 0.91986621 & 0.90939741 & 1 & 0.99992715 \cr
73 & 0.94733093 & 0.94345043 & 1.099525 & 1.1110782 & 1.0846523 & 1.0950853 \cr
79 & 1.1338934 & 1.1562237 & 0.90453403 & 0.8922209 & 1 & 0.99882039 \cr
83 & 1.0741723 & 1.0854551 & 0.8660254 & 0.84732408 & 1 & 0.99979996 \cr
89 & 0.84515425 & 0.82410673 & 0.87447463 & 0.85750248 & 0.89442719 & 0.88154899 \cr
97 & 1.0741723 & 1.0877289 & 0.92144268 & 0.90867892 & 0.8304548 & 0.80811684 \cr
101 & 0.98058068 & 0.97846254 & 0.94280904 & 0.93032086 & 1.0198039 & 1.0229108 \cr
103 & 1.1677484 & 1.1976448 & 0.87333376 & 0.855721 & 1 & 1.0004009 \cr
107 & 0.84515425 & 0.82186438 & 1.183216 & 1.2153554 & 1 & 1.0009282 \cr
109 & 0.91287093 & 0.89933354 & 1.1577675 & 1.1844329 & 0.94686415 & 0.94015124 \cr
113 & 0.92393644 & 0.9146531 & 0.9486833 & 0.93966595 & 1.1313708 & 1.1534106 \cr
127 & 0.93933644 & 0.93052596 & 0.98449518 & 0.98005032 & 1 & 0.99904006 \cr
131 & 1.1470787 & 1.171545 & 1.1208971 & 1.1413931 & 1 & 0.99916309 \cr
137 & 1.052079 & 1.0603352 & 1.0219806 & 1.0285831 & 1.1744404 & 1.2066518 \cr
139 & 0.93094934 & 0.91532106 & 1.0975994 & 1.1176423 & 1 & 1.0000469 \cr
149 & 1.069045 & 1.0833831 & 0.86855395 & 0.84844439 & 0.91064169 & 0.89706709 \cr
}}}
\botcaption{Table 1}
A table in support of Conjecture 2, comparing $R_p$ v.s. $R_p(T)$ for the three 
elliptic curves $E_{11}$, $E_{19}$, $E_{32}$
($T$ equal to 333605031, 263273979, 930584451 respectively).
More of this data, for $p<2000$, is depicted in the figures below.
The $0$ entries for $p=11$ and $p=19$ are explained by the fact that
we are restricting ourselves to twists with even functional equation,
$w_E \chi_d(-N)=1$. Hence for $E_{11}$ and $E_{19}$, we are
only looking at twists with $\chi_d(11) = \chi_d(19) = -1$.
\endcaption
\noindent
\endinsert
\bigskip

\midinsert
{
\centerline{\psfig{figure=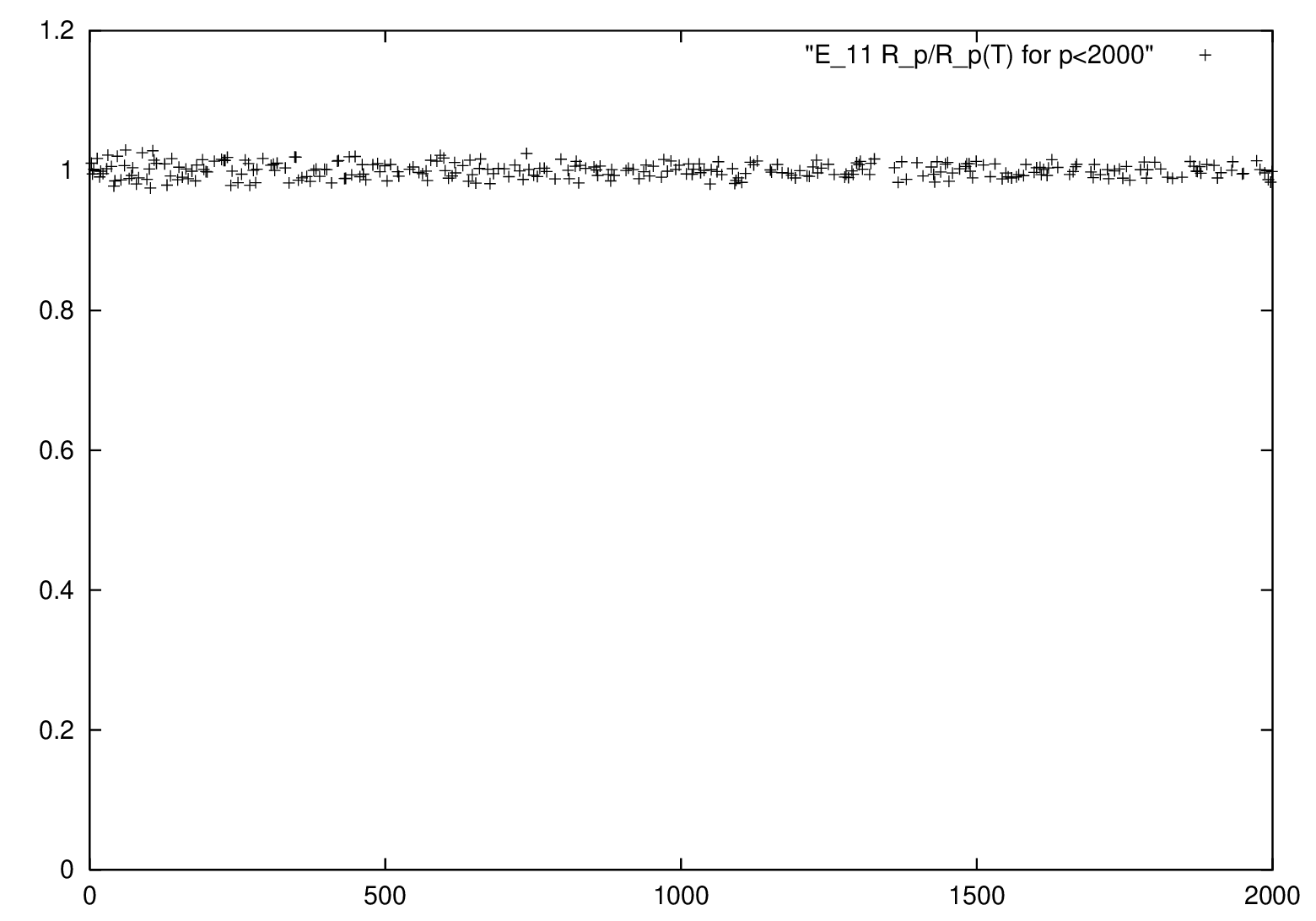,width=2.5in}
            \psfig{figure=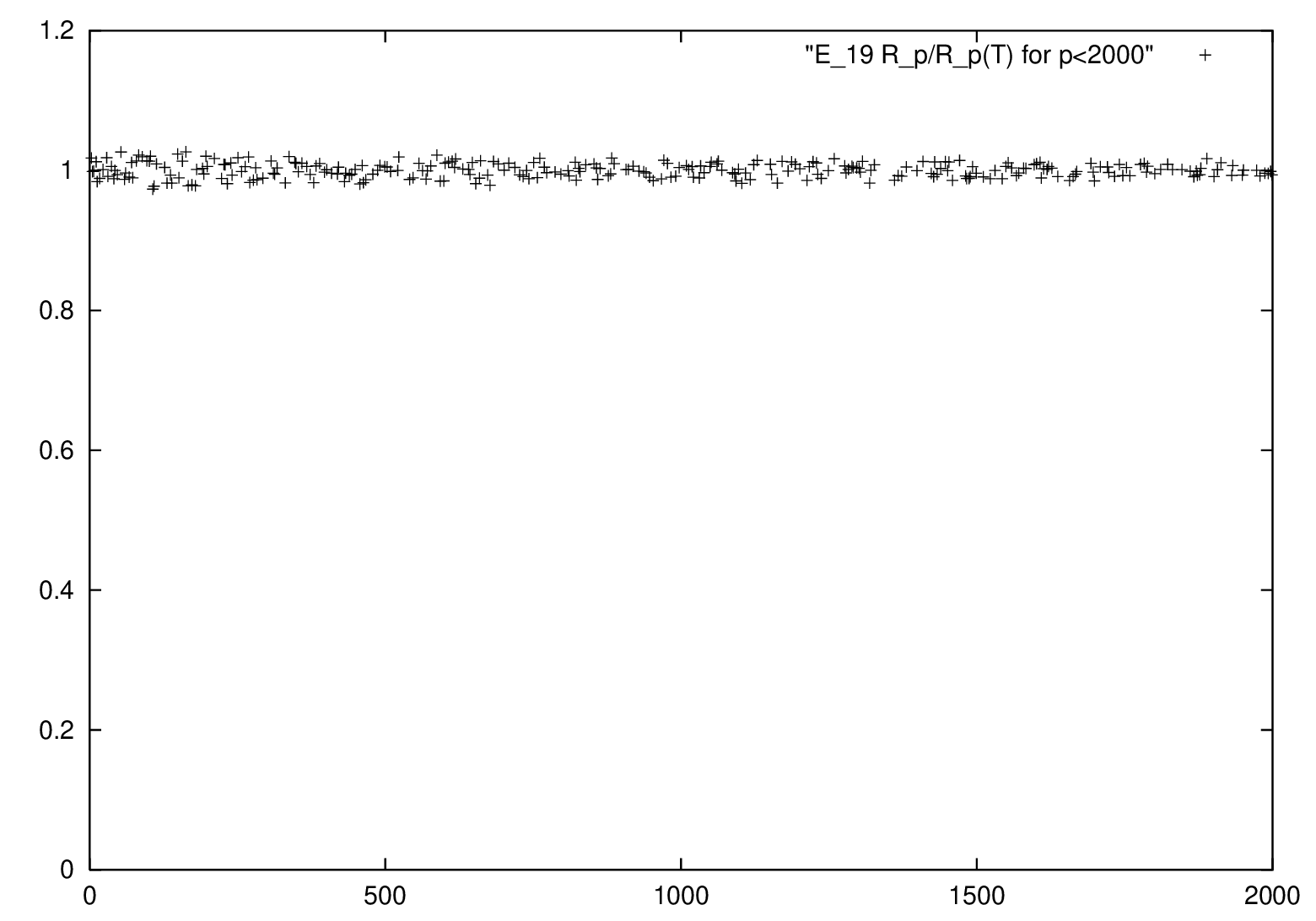,width=2.5in}}
\centerline{\psfig{figure=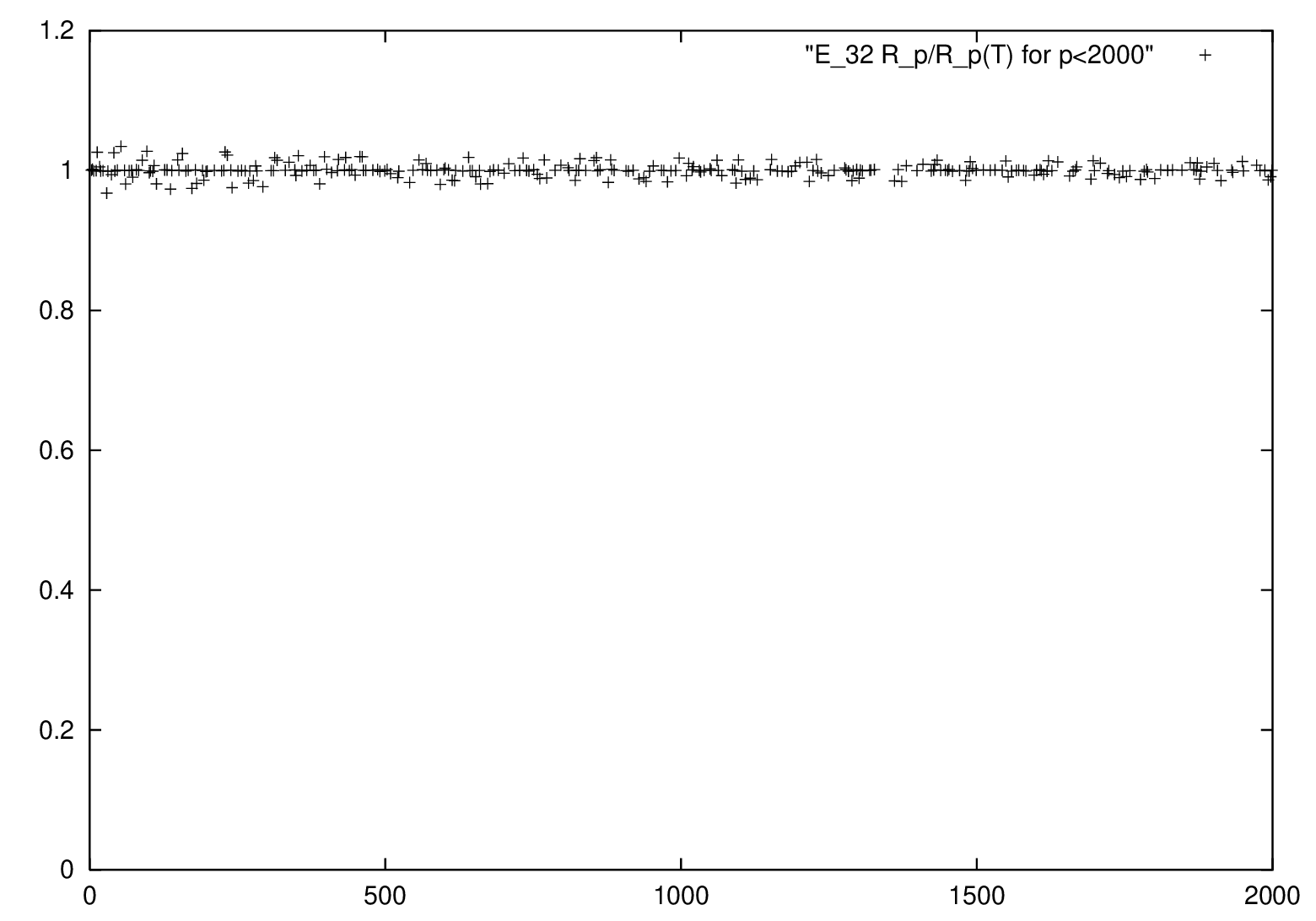,width=2.5in}}
}
\topcaption{Figure 1}
Pictures depicting $R_p/R_p(T)$, for $p<2000$, $T$ as in Table 1.
\endcaption
\endinsert

\midinsert
{
\centerline{\psfig{figure=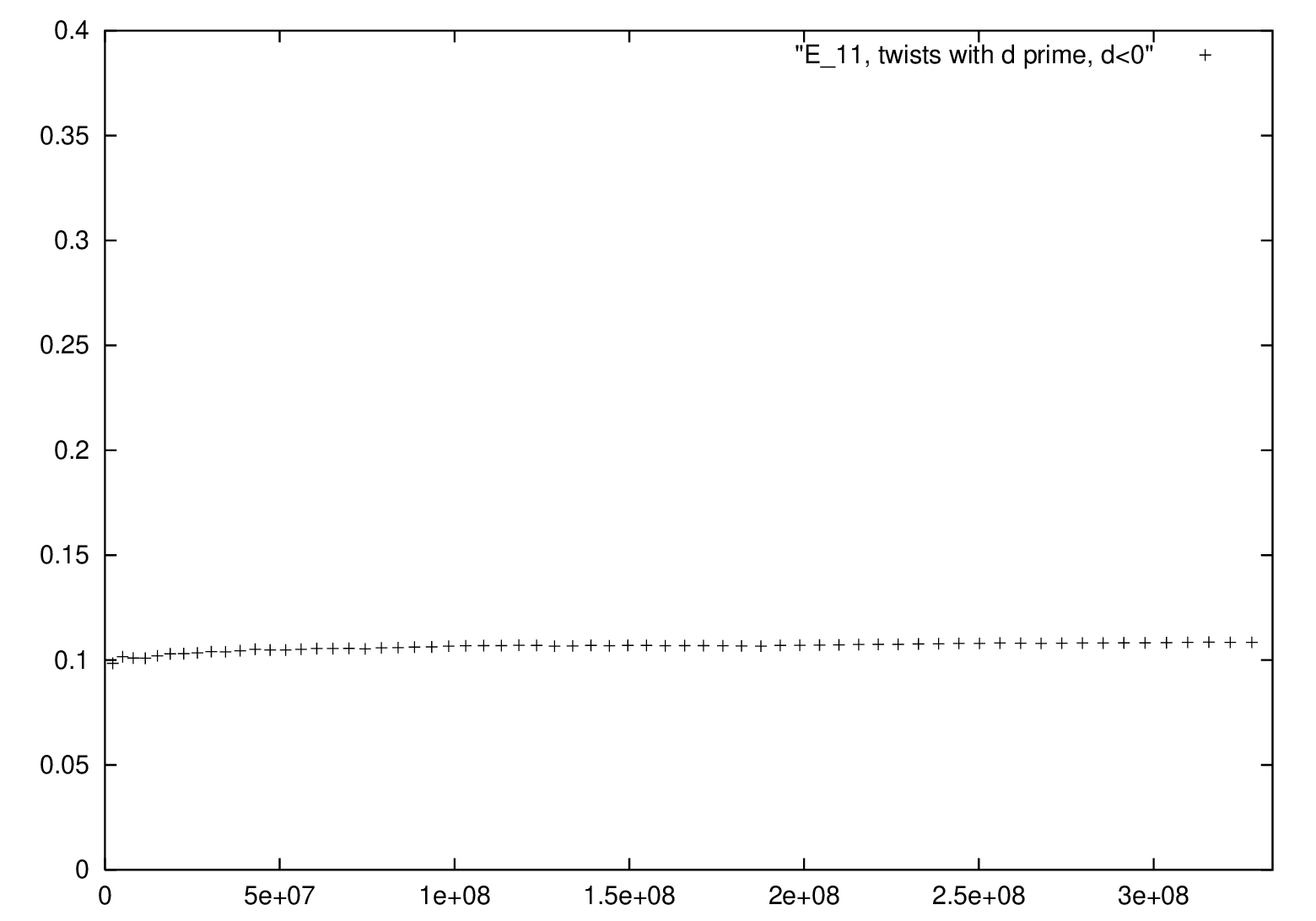,width=2.5in}
            \psfig{figure=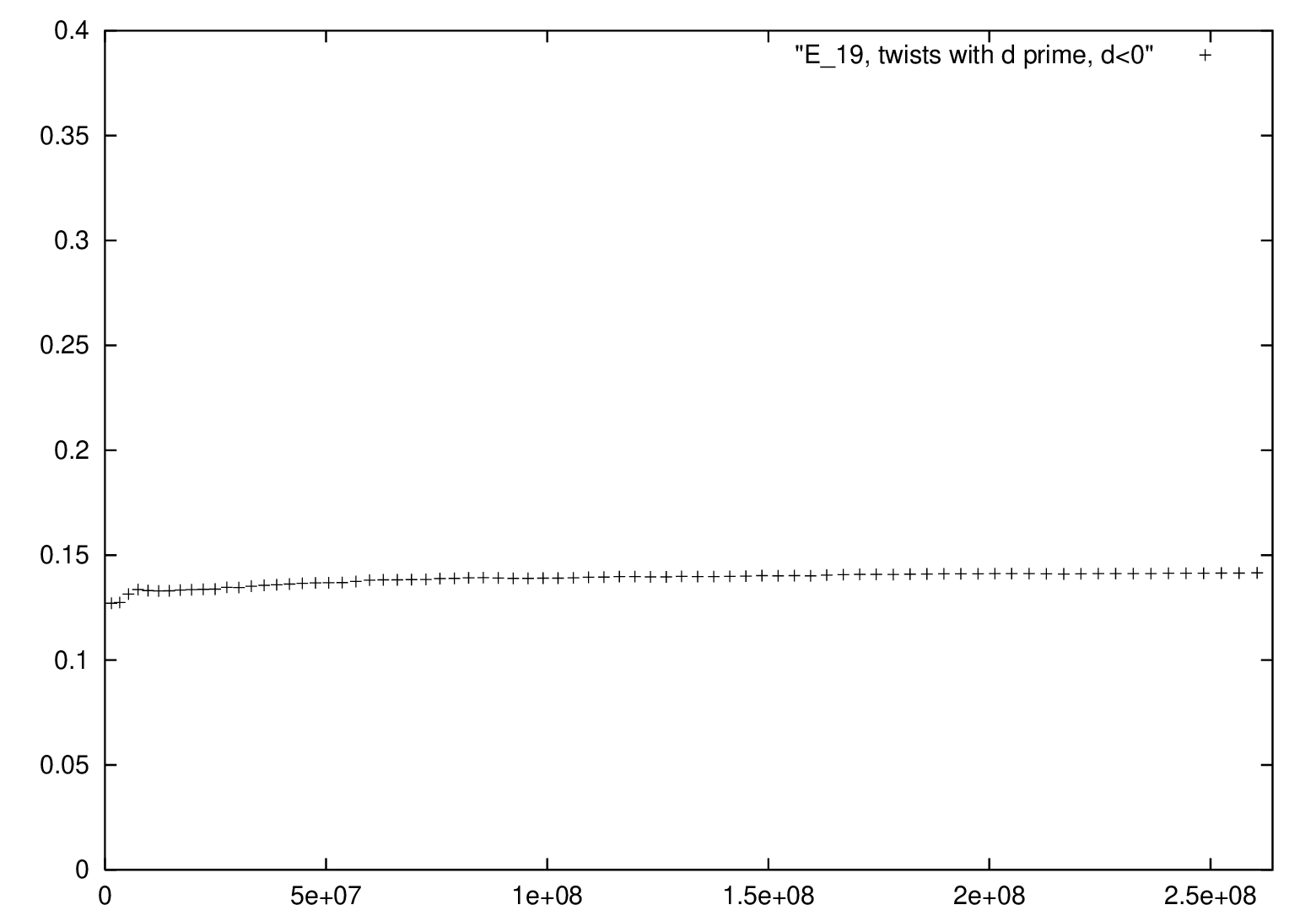,width=2.5in}}
\centerline{\psfig{figure=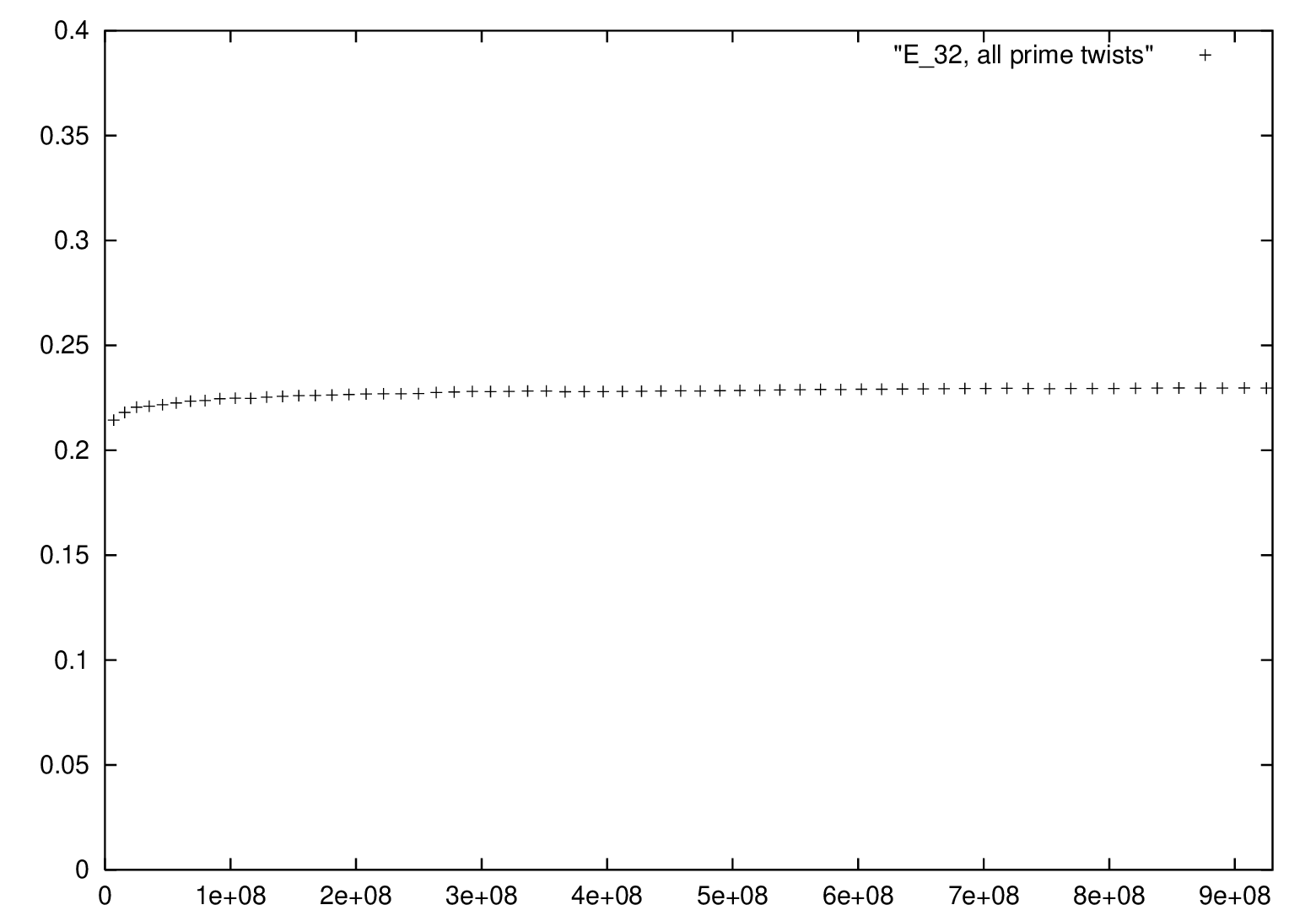,width=2.5in}
            \psfig{figure=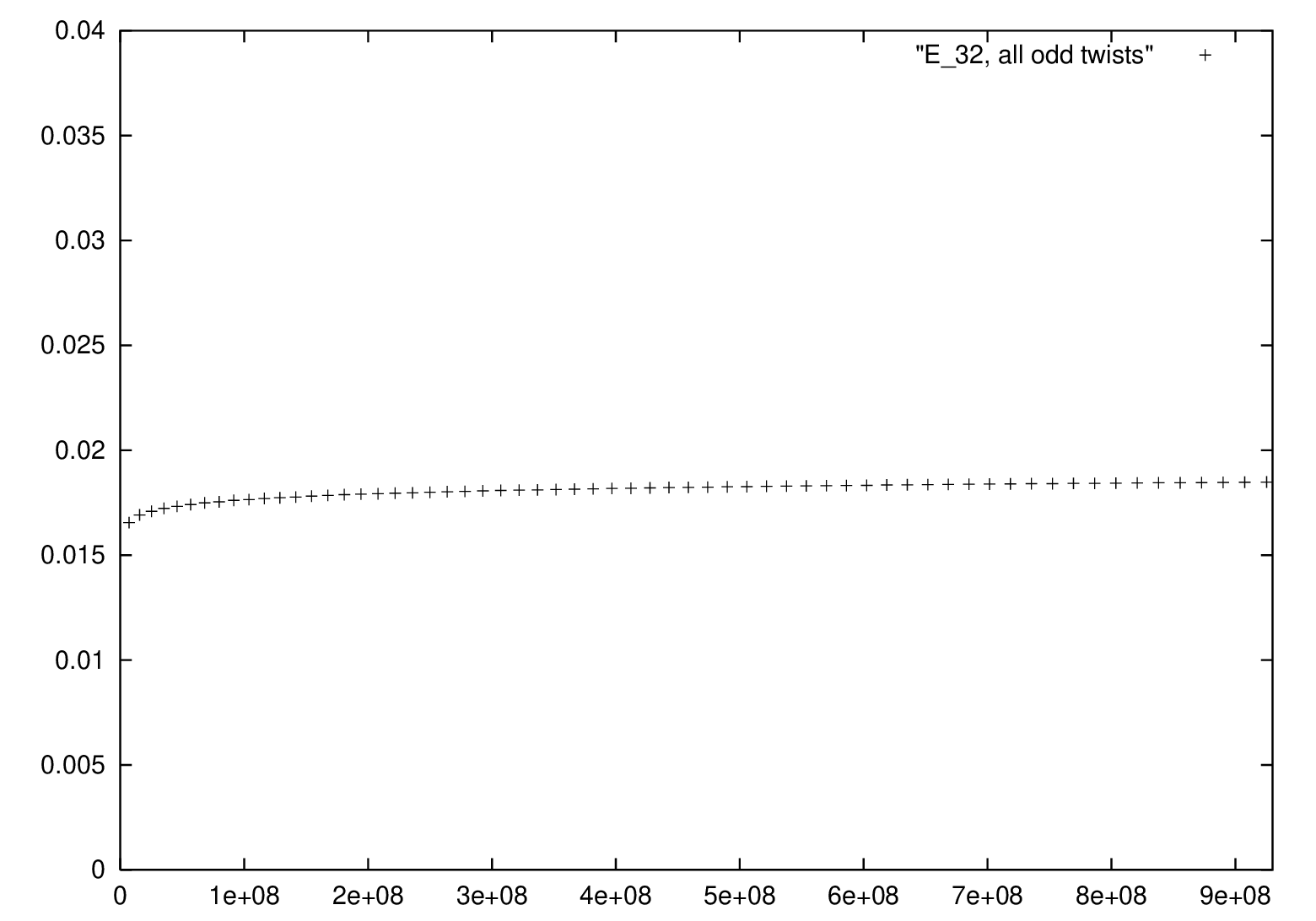,width=2.5in}}
}
\topcaption{Figure 2}
Figures in support of Conjecture 1. These depict the l.h.s. of (25)
divided by $T^{3/4}(\log T)^{-5/8}$. For the level
11 and 19 curves, we only looked at twists with $d<0$, $d$ prime, even
functional equation.
We also depict the l.h.s of (23) divided by $T^{3/4}(\log T)^{11/8}$
for the level 32 curve and odd $d$. While the pictures are reasonably flat, 
$\log(T)$ is almost constant for most of the interval in question.
The flatness we are observing reflects the main dependence on $T^{3/4}$.
\endcaption
\endinsert

\midinsert
{
\centerline{\psfig{figure=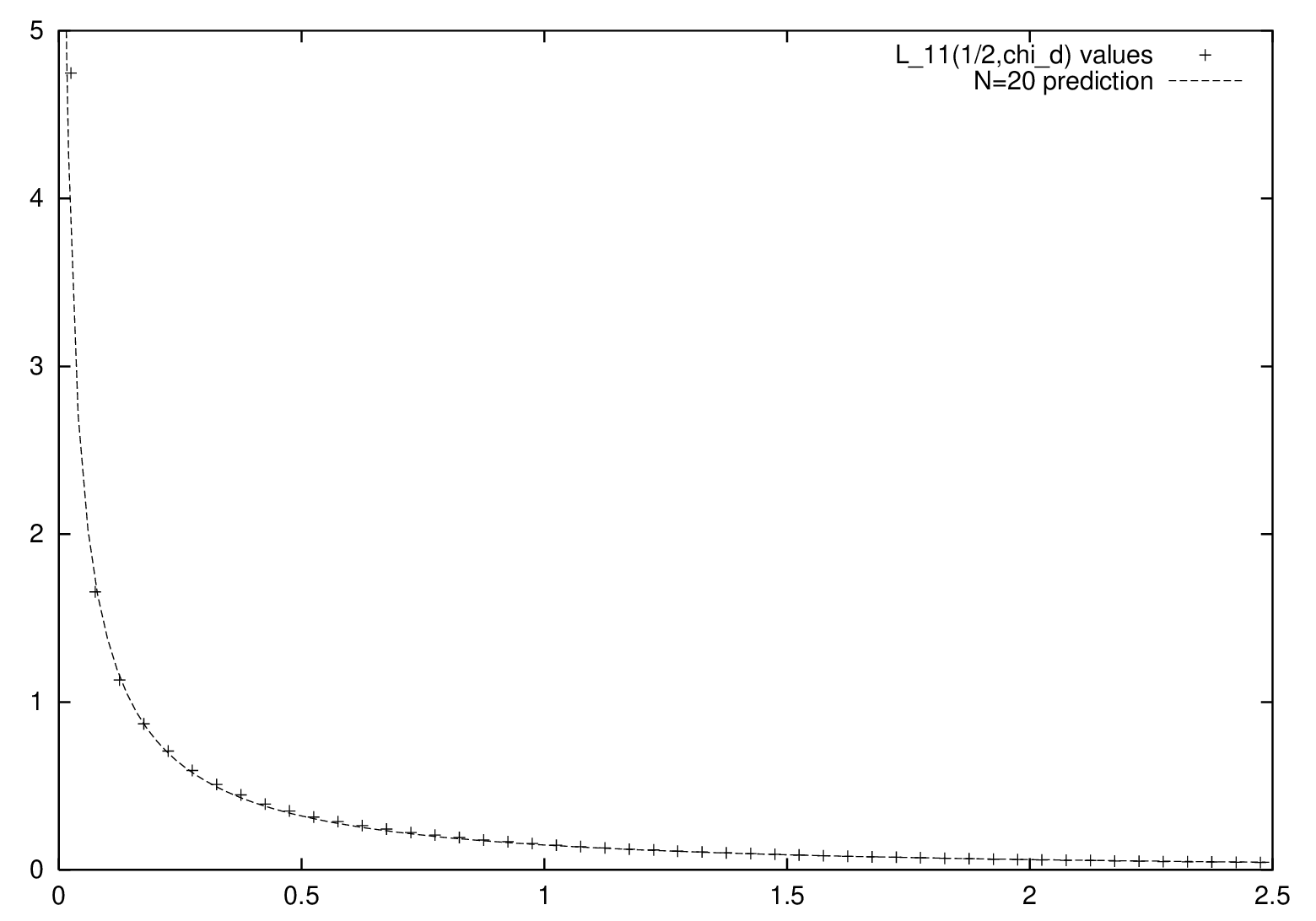,width=2.5in}
            \psfig{figure=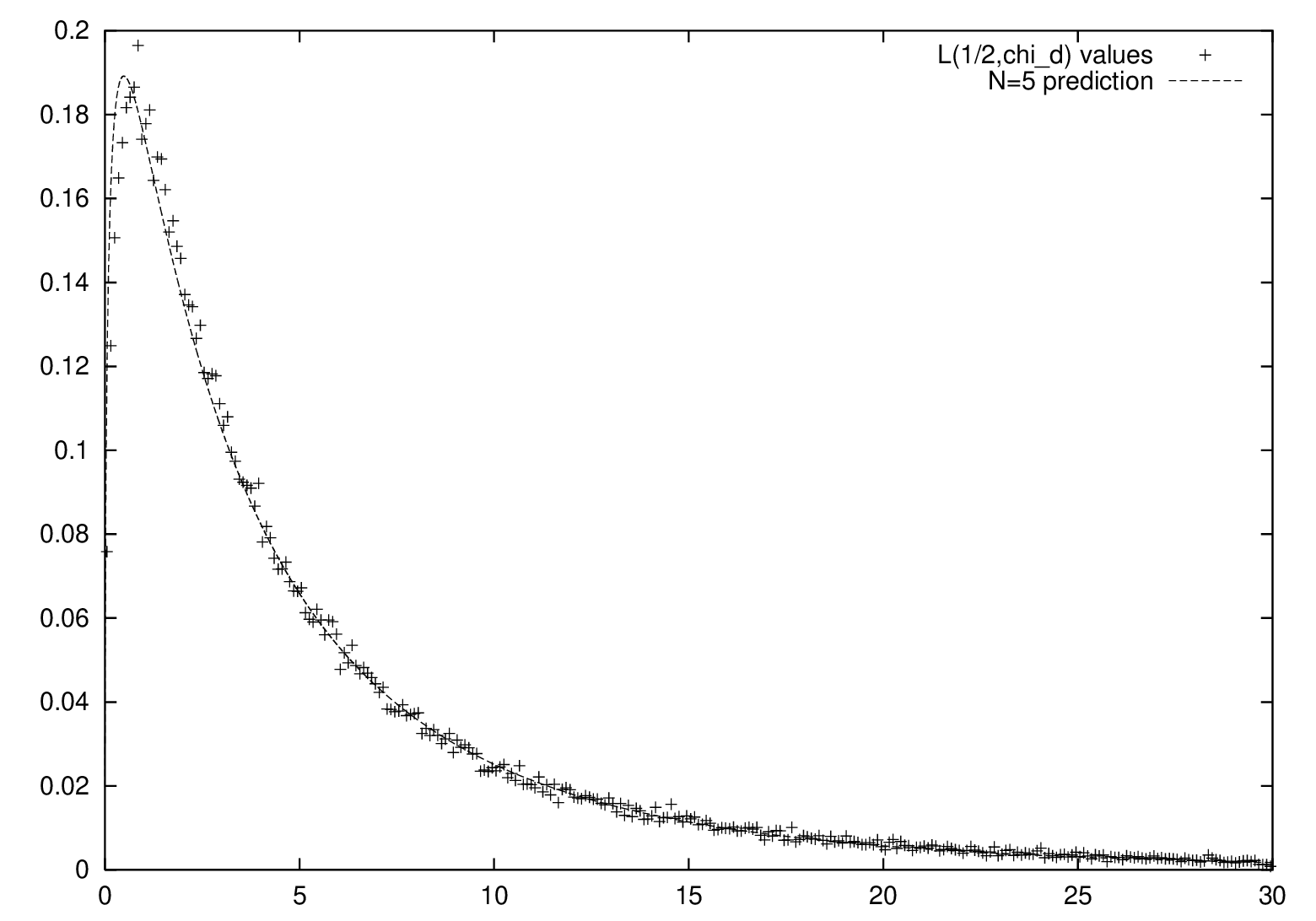,width=2.5in}}
}
\topcaption{Figure 3}
The first picture depicts the value distribution of $L_{E_{11}}(1/2,\chi_d)$,
for prime $|d|$,  $ -788299808 < d < 0 $, even functional equation, 
compared to $P_O(N,x)$, with 
$N=20$. For contrast, we depict, in the second picture, the value distribution 
of $L(1/2,\chi_d)$ (Dirichlet $L$-functions)
for all fundamental $800000 < |d| < 1000000$. Here, the Katz-Sarnak philosophy
predicts a Unitary Symplectic family, and so we compare the data against
$P_{USp}(N,x)$, $N=5$. In these pictures, we have renormalized the $L$ values
so as to have the same means as $P_O(20,x)$, and $P_{USp}(5,x)$ respectively,
and have not incorporated the $a_k$ values into the pictures. 
\endcaption
\endinsert

\enddocument